
\baselineskip=14pt
\parskip=10pt

 \font\sixrm=cmr6

\magnification=\magstephalf
\def\W{{\cal W}}

\def\1{{\overline{1}}}
\def\2{{\overline{2}}}
\parindent=0pt
\overfullrule=0in

\def\frac#1#2{{#1 \over #2}}

\centerline
{
\bf
 Doron Gepner's Statistics on Words in $\{1,2,3\}^{*}$ is (Most Probably) Asymptotically Logistic
}

\bigskip
\centerline
{\it  Doron ZEILBERGER}
\bigskip
\qquad 
{\it 
Dedicated to my friend and hero, Doron Gepner (b. March 31, 1956), on his $60^{th}$ birthday.
``Doron le Doron me-Doron''.
}
\bigskip

{\bf Preface}

I first met Doron Gepner in 1980, when he was a Physics graduate student at the Weizmann Institute of Science,
and I was a young {\it khoker bakhir}. Already then Doron was a legend, since he was the first person  in Israel,
as far as I know, to have solved {\it Rubik's cube} completely from scratch, using
group-theoretical methods. I was so impressed that I asked him to present a guest-lecture
in my graduate combinatorics class, and the students loved it.

Doron then went on to do seminal work in theoretical physics, 
that, unfortunately, is over my head.
But the part that is really interesting to me is his current work, greatly generalizing the celebrated {\it Rogers-Ramanujan} identities,
and giving lots of new insight. I am sure that this work will lead to many future gems.

The purpose of the present note is to present a {\it present}
({\it Doron} in Hebrew)\footnote{${}^1$}{\sixrm Exactly thirty years ago, on March 31, 1986, my wife Jane and I were at Doron Gepner's
30th birthday party, (in Princeton) that Ida (Doron's wife) organized in their place. Another guest was
a colleague of Doron, an Egyptian postdoc, and we pointed out to him
that we have the same name, and that
it means a ``gift'', (presumably ``God's gift''), to which he
retorted  ``  in your cases it seems to be the devil's gifts''.},
from one Doron to another, by paying an old {\it debt}. In 1987, when he was
a postdoc at Princeton University, Doron introduced a new {\it permutation statistic} (see below),
that came up in his work in string theory and conformal field theory.
It so happened that at the  time,
my friend and collaborator, the eminent French combinatorialist, Dominique Foata, visited me.
Foata is the world's greatest expert on permutation (and word-) statistics 
(and coined the term!), so it was only natural that
we both got intrigued and tried to investigate Gepner's new statistics,
that we christened {\it gep}, in analogy with the classical statistics {\it inv}, {\it maj}, and {\it des} (see below).
We had some preliminary results, but not enough for a paper.
This was due to the fact that my beloved servant, Shalosh B. Ekhad, was not yet born. Now, almost thirty years later, it
is a good opportunity to revisit Doron Gepner's difficult statistics and harness the full power of
my silicon servant, and of Maple, to study it seriously.

{\bf Important note}:
All the results in this paper were gotten by using the Maple package {\tt GEPNER.txt}, available, free of charge, from
the url

{\tt http://www.math.rutgers.edu/\~{}zeilberg/tokhniot/GEPNER.txt} \quad .

Sample input and output files may be gotten from the front of this article:

{\tt http://www.math.rutgers.edu/\~{}zeilberg/mamarim/mamarimhtml/gepner.html} \quad .

{\bf A crash course on Permutation and Word Statistics}

A permutation statistics is an integer-valued function on the set of permutations. The most famous one is the
{\bf number of inversions}, $inv(\pi)$, (that shows up in the definition of the determinant of a square matrix)
$$
inv(\, \pi_1 \dots \pi_n \, ):=\sum_{1 \leq i<j \leq n}   \chi (\pi_i>\pi_j) \quad ,
$$
where $\chi(S)$ is $1$ or $0$, according to whether $S$ is true or false, respectively.

Almost as famous is Major Percy Alexander MacMahon's statistics, $maj(\pi)$, called the ``major index''
$$
maj(\, \pi_1 \dots \pi_n \, ):=\sum_{i=1}^{n-1} \, i \, \chi (\pi_i>\pi_{i+1}) \quad .
$$
The generating functions according to $inv$ and $maj$ are both given by 
$[n]!:=1 (1+q)(1+q+q^2) \cdots (1+q+\dots +q^{n-1})=(1-q)\dots (1-q^n)/(1-q)^n)$ 
as proved by Netto and MacMahon respectively.
In particular the permutation statistics $inv$ and $maj$ are {\it equally distributed}.

Permutations of length $n$ may be viewed as {\it words} in the ``alphabet'' $\{1,2, \dots, n\}$ with
exactly one occurrence of each letter. The above definitions of $inv$ and $maj$ make perfect sense
when defined on words, of any length, in the same alphabet, where repetitions (and omissions) are welcome.

Let $\W(a_1, \dots, a_n)$ be the set of words in the alphabet $\{1, \dots , n\}$ 
with $a_1$ occurrences of $1$, $a_2$ occurrences of $2$, $\dots$, $a_n$ occurrences  of $n$.
MacMahon proved (Theorems 3.6 and 3.7 in [A])
$$
\sum_{w \in \W(a_1, \dots, a_n)} q^{inv(w)} = \frac{[a_1 + \dots + a_n]!}{[a_1]! \dots [a_n]!} \quad,
$$
$$
\sum_{w \in \W(a_1, \dots, a_n)} q^{maj(w)} = \frac{[a_1 + \dots + a_n]!}{[a_1]! \dots [a_n]!} \quad,
$$
(where, as mentioned above, $[m]!:=(1-q)(1-q^2) \cdots (1-q^m)/(1-q)^m$).

In particular they are still {\it equally distributed}, and Dominique Foata ([F]) gave a gorgeous 
bijective proof.

\vfill\eject

{\bf Asymptotic Normality}

Most (but not all) combinatorial statistics (naturally parametrized by one or more integer parameters), for example
tossing a (fair or loaded) coin $n$ times and observing the number of Heads, are {\it asymptotically normal}, that
means that (for the sake of simplicity let's only consider the one parameter case), if you call the sequence $X_n$,
figure out its {\it mean}, $\mu_n$, (usually extremely easy) (aka as average, aka as expectation), call it $\mu_n$, and its
{\it variance}, $m_2(n)$ (also, usually, fairly easy), and define the {\it standardized} sequence of random variables
$$
Z_n :=\frac{X_n -\mu_n}{\sqrt{m_2(n)}} \quad,
$$
then the sequence $\{Z_n\}$, converges, in distribution, to the good-old normal distribution 
(aka Gaussian distribution) whose probability density function is, famously,
$e^{-x^2/2}/\sqrt{2\pi}$. A good way to prove this is to discover explicit expressions for the moments (about the mean),
$m_r(n)$ (or at least the leading terms), and prove that the {\it standardized moments}, $m_r(n)/m_2(n)^{r/2}$ tend, as $n$ goes to infinity,
to the moments of the standard normal distribution, that equal $0$ when $r$ is odd and $1 \cdot 3 \cdots (r-1)$, when $r$ is even.
This approach can be often taught to a computer, see [Z1][Z2]. 
The asymptotic normality of $inv$, for the two-lettered case,
was first proved  by Mann and Whitney [MW], and for the  general case by Persi Diaconis[D] (and reproved in [CJZ]).

{\bf Enter Doron Gepner's Statistics}

One way to look at an inversion in a word $w_1 \dots w_m$ is as the number of {\bf pairs}
of letters $w_i w_j$ with $1 \leq i<j \leq m$,
whose {\it reduction} is the permutation of length $2$, $21$, the only {\bf odd} permutation of length $2$.

This leads naturally to the analog for three-letter subwords. Given a word 
$w=w_1 \dots w_m$ (in any alphabet)
consider the set of {\bf triples} $1 \leq i<j < k \leq m$ such that $w_iw_jw_k$ reduces to one of
the three odd permutations of length $3$, namely one of the members of $\{132,213,321\}$.

And indeed, this came up (naturally!), in Doron Gepner's[G] deep work in conformal field theory, and leads to
the following definition.

{\bf Definition.} The Gepner statistics on a word $w=w_1 \dots w_m$ on any finite, totally ordered, alphabet,  
denoted by $gep(w)$, is defined by
$$
gep(w_1 \dots w_m):=
\sum_{1 \leq i<j <k \leq m} \chi (w_i<w_k<w_j \quad OR \quad w_j<w_i<w_k  \quad OR \quad w_k<w_j<w_i) \quad .
$$

It seems that the Gepner statistics, $gep$, is much harder to study than the inversion  number,$inv$.  and it is
extremely unlikely that a `nice' (or even `ugly') explicit formula for the generating function exists.
But it is still interesting to find out whether it is asymptotically normal, and if not, to determine
the limiting distribution.

Failing an explicit formula for the generating function, 
the best that we can hope for is an {\it algorithm} to compute the first
few terms of the generating functions for the Gepner statistics. Using the method of [Z1] and [Z2], we
can then infer explicit (rigorously proved!) {\bf polynomial} expressions for the first few moments, and
try to see what is going on. It is easy to see by a (fully rigorous) {\it handwaving argument} that  each moment
is always {\it some} polynomial, and it is also easy to bound the degree of the $r$-th moment to
be $3r$  (in the cases considered here). The degree of the polynomial expression, in $n$, of
the $r$-th moment-about-the mean is lower, in fact, in the two cases considered below,
it happens to be $2r$. Since a polynomial of degree $d$ is uniquely determined by $d+1$ distinct values, an experimental-yet-rigorous
approach would be to find numerical values and then ``fit the data''.

Due to the complexity of the Gepner statistics (and our finite time) we will only  consider the cases
of permutations, i.e. $\W(1, \dots, 1)$ (with $n$ $1$s), and
words in the three-letter alphabet $\{1,2,3\}$ with $n$ occurrences of each of the
letters $1$, $2$, and $3$, i.e. $\W(n,n,n)$.

{\bf First Surprise: The Gepner Statistics on Permutations is NOT Asymptotically Normal}

Define the {\it Gepner polynomials of the first kind}, $G_n(q)$, to be the polynomial
$$
G_n(q):=\sum_{\pi \in S_n} q^{gep(\pi)} \quad .
$$
Here are the first $8$ Gepner polynomials
$$
G_1(q)=1 \quad , \quad
G_2(q)=2 \quad , \quad
G_3(q)=3q+3 \quad,\quad
G_4(q)=4q^4+16q^2+4 \quad,
$$
$$
G_5(q)=  5\,{q}^{10}+25\,{q}^{7}+25\,{q}^{6}+10\,{q}^{5}+25\,{q}^{4}+25\,{q}^{3}+5 \quad,
$$
$$
G_6(q)=  6\,{q}^{20}+36\,{q}^{16}+72\,{q}^{14}+138\,{q}^{12}+216\,{q}^{10}+138\,{q}^{8}+72\,{q}^{6}+36\,{q}^{4}+6 \quad ,
$$
$$
G_7(q)=  7\,{q}^{35}+49\,{q}^{30}+98\,{q}^{27}+49\,{q}^{26}+98\,{q}^{25}+98\,{q}^{24}+147\,{q}^{23}+441\,{q}^{22}+308\,{q}^{21}+196\,{q}^{20}
$$
$$
+490\,{q}^{19}  +539\,{q}^{18}+539\,{q}^{17}+490\,{q}^{16}+196\,{q}^{15}+308\,{q}^{14}+441\,{q}^{13}+147\,{q}^{12}+98\,{q}^{11}+98\,{q}^{10}
$$
$$
+49\,{q}^{9}+98\,{q}^{8}+49\,{q}^{5}+7 \quad,
$$
$$
G_8(q)= 
8\,{q}^{56}+64\,{q}^{50}+128\,{q}^{46}+288\,{q}^{44}+128\,{q}^{42}+896\,{q}^{40}+1344\,{q}^{38}+1600\,{q}^{36}+3200\,{q}^{34}
$$
$$
+4792\,{q}^{32}+4352\,{q}^{30}+ 6720\,{q}^{28}+4352\,{q}^{26}+4792\,{q}^{24}+3200\,{q}^{22}+1600\,{q}^{20}+1344\,{q}^{18}+896\,{q}^{16}
$$
$$
+128\,{q}^{14}+288\,{q}^{12}+128\,{q}^{10}+64\,{q}^{6}+8
\quad.
$$

By looking at the Taylor expansion around $q=1$, one can get the so-called factorial moments, and from them the usual moments,
and from them the moments about the mean (see [Z1]), and Shalosh B. Ekhad found  that the average is
$$
\mu_1(n) =\frac{1}{12} \,n \left( n-1 \right)  \left( n-2 \right)  \quad ,
$$
(of course, this is obvious, since it must be ${{n} \choose {3}}/2$ [why?]).
Less trivially (but still humanly doable), it found that the variance, $m_2(n)$, is given by the
following polynomial expression
$$
m_2(n)={\frac {1}{72}}\,{n}^{2} \left( n-1 \right)  \left( n-2 \right)  \quad,
$$
in particular, the standard deviation is $O(n^2)$, while the average is $O(n^3)$, so we have, at least, {\it concentration about the mean}.

Of course, by symmetry, $m_3(n)=0$, but much less trivially, we have the following {\bf rigorously proved} formula
for the fourth moment
$$
m_4(n)={\frac {1}{43200}}\,{n}^{3} \left( n-1 \right)  \left( n-2 \right)  \left( 29\,{n}^{3}-71\,{n}^{2}-26\,n-16 \right)  \quad.
$$
From this it follows that the {\it standardized fourth moment}, aka as the {\bf kurtosis}, equals
$$
\alpha_4(n)=
\frac{
{\frac {1}{43200}}\,{n}^{3} \left( n-1 \right)  \left( n-2 \right)  \left( 29\,{n}^{3}-71\,{n}^{2}-26\,n-16 \right)
}
{
\left ({\frac {1}{72}}\,{n}^{2} \left( n-1 \right)  \left( n-2 \right)  \right )^2 \quad
}
$$
and taking the limit as $n \rightarrow \infty$, we get that the {\bf limiting kurtosis} is
$$
\frac{87}{25} = 3.48 \quad,
$$
and {\bf not} $3$, so the Gepner statistics is {\bf not} asymptotically normal,
and happens to be {\it leptokurtic} (statisticians' big word for saying that the kurtosis exceeds $3$).
We have no clue about the limiting distribution, and would love to know it.

The difficulty, from a computational point of view, is that we have no efficient way to get more
Gepner polynomials than by very naive  brute force. My former PhD student, Brian Nakamura([N]),
found more efficient algorithms for related problems, but this needs more work.
On the  other hand, for the case below, where we look at the Gepner statistics defined
on $\W(n,n,n)$, we can do much better than brute force, and squeeze out explicit expressions for more moments.

{\bf Second Surprise: The Gepner Statistics on Words in $\{1,2,3\}^{*}$ is Asymptotically Logistic}

Define the {\it Gepner polynomials of the second kind}, $g_n(q)$, to be the polynomial
$$
g_n(q):=\sum_{w \in \W(n,n,n) } q^{gep(w)} \quad .
$$
(Recall that $\W(n,n,n)$ is the set of words with $n$ occurrences each of the letters $1$, $2$, and $3$).

Here are the first five Gepner  polynomials of the second kind.

$$
g_1(q)=3q+3 \quad, \quad g_2(q)=6\,{q}^{8}+21\,{q}^{6}+36\,{q}^{4}+21\,{q}^{2}+6 \quad ,
$$
$$
g_3(q)=9\,{q}^{27}+27\,{q}^{24}+108\,{q}^{21}+264\,{q}^{18}+432\,{q}^{15}+432\,{q}^{12}+264\,{q}^{9}+108\,{q}^{6}+27\,{q}^{3}+9
\quad ,
$$
$$
g_4(q)=
12\,{q}^{64}+36\,{q}^{60}+108\,{q}^{56}+336\,{q}^{52}+870\,{q}^{48}+2016\,{q}^{44}+4041\,{q}^{40}+6252\,{q}^{36}+7308\,{q}^{32}
\quad ,
$$
$$
+6252\,{q}^{28}+4041\,{q}^{24
}+2016\,{q}^{20}+870\,{q}^{16}+336\,{q}^{12}+108\,{q}^{8}+36\,{q}^{4}+12 \quad ,
$$
$$
g_5(q)=
15\,{q}^{125}+45\,{q}^{120}+135\,{q}^{115}+330\,{q}^{110}+900\,{q}^{105}+2115\,{q}^{100}+4710\,{q}^{95}+10230\,{q}^{90}+21195\,{q}^{85}+40290\,{q}^{80}
$$
$$
+69423\,{q}^{75}+102780\,{q}^{70}+126210\,{q}^{65}+126210\,{q}^{60}+102780\,{q}^{55}+69423\,{q}^{50}+40290\,{q}^{45}+21195\,{q}^{40}+10230\,{q}^{35}
$$
$$
+4710\,{q}^{30}+2115\,{q}^{25}+900\,{q}^{20}+330\,{q}^{15}+135\,{q}^{10}+45\,{q}^{5}+15 \quad .
$$

{\bf How to get More Gepner Polynomials?}

Consider, for motivation, computing the generating function, according to $inv$, of the words in
$\W(a_1,a_2,a_3)$, i.e. how can we generate many terms of the polynomials
$$
f(a_1,a_2,a_3):=\sum_{w \in \W(a_1,a_2,a_3)} q^{inv(w)} \quad .
$$

Suppose that we did not know MacMahon's explicit formula $[a_1+a_2+a_3]!/([a_1]![a_2]![a_3]!)$ for it.
A natural approach would be to use {\it Dynamical Programming} and try to express $f(a_1,a_2,a_3)$ in terms
of $f$ with smaller arguments.

Note that for $w' \in \W(a_1-1,a_2,a_3)$,
$$
inv(w'1)=inv(w')+ a_2+a_3 \quad ,
$$
since appending a $1$ to $w'$ creates $a_2$ new inversions due to the $a_2$ $2$s, and $a_3$ new inversions due to
the $a_3$ $3$'s.

Similarly, for $w' \in \W(a_1,a_2-1,a_3)$,
$$
inv(w'2)=inv(w')+ a_3 \quad ,
$$
since appending a $2$ to $w'$ creates  $a_3$ new inversions due to the $a_3$ $3$'s.

Finally, for $w' \in \W(a_1,a_2,a_3-1)$,
$$
inv(w'3)=inv(w') \quad ,
$$
since appending a $3$ to $w'$ does not create any new inversions.
Since
$$
\W(a_1,a_2,a_3)= \W(a_1-1,a_2,a_3) 1 \cup \W(a_1,a_2-1,a_3) 2 \cup \W(a_1,a_2,a_3-1) 3 \quad ,
$$
\vfill\eject
we have
$$
f(a_1,a_2,a_3):=\sum_{w \in \W(a_1,a_2,a_3)} q^{inv(w)}
$$
$$
=\sum_{w' \in \W(a_1-1,a_2,a_3)} q^{inv(w')+a_2+a_3} 
+\sum_{w' \in \W(a_1,a_2-1,a_3)} q^{inv(w')+a_3} 
\sum_{w' \in \W(a_1,a_2,a_3-1)} q^{inv(w')} 
$$
$$
q^{a_2+a_3}f(a_1-1,a_2,a_3)+ q^{a_3}f(a_1,a_2-1,a_3)+ f(a_1,a_2,a_3-1) \quad .
$$
We have just established the recurrence
$$
f(a_1,a_2,a_3)=q^{a_2+a_3}f(a_1-1,a_2,a_3)+ q^{a_3}f(a_1,a_2-1,a_3)+ f(a_1,a_2,a_3-1) \quad ,
$$
that would enable us to crank out many $f(n,n,n)$.

In fact, this recurrence can be also used to prove MacMahon's formula, by verifying that $[a_1+a_2+a_3]!/([a_1]![a_2]![a_3]!)$
also satisfies the same recurrence (a trivial high-school algebra verification, even for humans), and that the boundary conditions
agree.

How can we generalize this argument for computing
$$
F(a_1,a_2,a_3):=\sum_{w \in \W(a_1,a_2,a_3)} q^{gep(w)} \quad .
$$

Let's use, once again the decomposition
$$
\W(a_1,a_2,a_3)= \W(a_1-1,a_2,a_3) 1 \cup \W(a_1,a_2-1,a_3) 2 \cup \W(a_1,a_2,a_3-1) 3 \quad .
$$

Let's consider a member, $w$, of $\W(a_1,a_2,a_3)$, that ends with the letter $1$, so
$w=w'1$ for $w' \in \W(a_1-1,a_2,a_3)$. Appending $1$ to $w'$ introduces new contributions to the $gep$ statistics
due to all the occurrences of the pairs $32$.

Let's consider a member, $w$, of $\W(a_1,a_2,a_3)$, that ends with the letter $2$, so
$w=w'2$ for $w' \in \W(a_1,a_2-1,a_3)$. Appending $2$ to $w'$ introduces new contributions to the $gep$ statistics
due to the occurrences of the pairs $13$.

Let's consider a member, $w$, of $\W(a_1,a_2,a_3)$, that ends with the letter $3$, so
$w=w'3$ for $w' \in \W(a_1,a_2,a_3-1)$. Appending $3$ to $w'$ introduces new contributions to the $gep$ statistics
due to the occurrences of the pairs $21$.

This forces us to keep track of the number of occurrences of the pairs $32$, $13$, and $21$.

So we have 
$$
gep(w'1)-gep(w')= \#32(w') \quad , \quad
gep(w'2)-gep(w')= \#13(w') \quad , \quad
gep(w'3)-gep(w')= \#21(w') \quad .
$$
We also have
$$
\#32(w'1)-\#32(w')=0 \quad, \quad
\#32(w'2)-\#32(w')=a_3 \quad, \quad
\#32(w'3)-\#32(w')=0 \quad, 
$$
$$
\#13(w'1)-\#13(w')=0 \quad, \quad
\#13(w'2)-\#13(w')=0 \quad, \quad
\#13(w'3)-\#13(w')=a_1 \quad, 
$$
$$
\#21(w'1)-\#21(w')=a_2 \quad, \quad
\#21(w'2)-\#21(w')=0 \quad, \quad
\#21(w'3)-\#21(w')=0 \quad.
$$

Defining the family of polynomials in $q$ and the three {\bf catalytic variables} $t_{32}$, $t_{13}$, and $t_{21}$:
$$
{\bf F} (a_1,a_2,a_3;q, t_{32},t_{13}, t_{21}):=\sum_{w \in \W(a_1,a_2,a_3)} q^{gep(w)} {t_{32}}^{\#32(w)}
{t_{13}}^{\#13(w)}{t_{21}}^{\#21(w)}
\quad ,
$$
we get the {\bf functional-recurrence} equation
$$
{\bf F} (a_1,a_2,a_3;q, t_{32},t_{13}, t_{21})=
$$
$$
t_{21}^{a_2} \cdot {\bf F} (a_1-1,a_2,a_3;q, qt_{32},t_{13}, t_{21})+
t_{32}^{a_3} \cdot {\bf F} (a_1,a_2-1,a_3;q, t_{32},qt_{13}, t_{21})+
t_{13}^{a_1} \cdot {\bf F} (a_1,a_2,a_3-1;q, t_{32},t_{13}, qt_{21}) \quad .
$$
At the {\it end of the day}, we plug-in $t_{32}=1$, $t_{13}=1$, $t_{21}=1$, and get our object of desire
$$
F(a_1,a_2,a_3;q)={\bf F} (a_1,a_2,a_3;q, 1,1, 1) \quad .
$$
and finally, we can get quite a few Gepner polynomials of the second kind from
$$
g_n(q)=F(n,n,n;q) \quad .
$$

But we can do better! If we are only interested in, say, the first $20$ moments, then we don't need the full polynomials.
Since the first $r$ factorial (and hence true) moments can be gotten from the first $r$ Taylor coefficients, about $q=1$, of $g_n(q)$,
we can do a change variables
$$
q=1+p \quad, \quad
t_{32}=1+x \quad , \quad
t_{13}=1+y \quad , \quad
t_{21}=1+z \quad , \quad
$$
and define
$$
{\bf H}(a_1,a_2,a_3; p,x,y,z):= {\bf F} (a_1,a_2,a_3;1+p, 1+x,1+y, 1+z) \quad .
$$
We get a corresponding functional recurrence equation for {\bf H}, and {\bf now} at each step, we can truncate,
and only retain the terms in $p,x,y,z$ of total degree $\leq r$, and of course, at the end of the day,
plug-in $x=0,y=0,z=0$. This enabled us to get many more truncated Gepner polynomials, and
enabled us to find explicit expressions for the first $12$ moments. This leads to

{\bf Theorem}: The average of the Gepner statistics, defined on the set of words in $\{1,2,3\}$ with
$n$ occurrences each of $1$, $2$, and $3$, is (of course)
$$
\mu_n=\frac{n^3}{2}    \quad .
$$
The variance is
$$
m_2(n)=\frac{n^4}{4} \quad.
$$
Of course, all the odd moments are zero. We also have, for the fourth through $12$th moments we have the following polynomial expressions.
$$
m_4(n)={\frac {1}{80}}\,{n}^{7} \left( -16+21\,n \right) \quad .
$$
$$
m_6(n)= {\frac {1}{448}}\,{n}^{9} \left( 279\,{n}^{3}-656\,{n}^{2}+512\,n-128 \right) \quad .
$$
$$
m_8(n)= {\frac {1}{1280}}\,{n}^{11} \left( 3429\,{n}^{5}-16480\,{n}^{4}+32512\,{n}^{3}-32512\,{n}^{2}+16384\,n-3328 \right)  \quad .
$$
$$
m_{10}(n)=
{\frac {1}{11264}}\,{n}^{13} \left( 3\,n-4 \right)   \cdot
$$
$$
\left( 68985\,{n}^{6}-469716\,{n}^{5
}+1391760\,{n}^{4}-2251584\,{n}^{3}+2072832\,{n}^{2}-1022208\,n+209920 \right) \quad .
$$
$$
m_{12}(n)=
{n}^{15} ( {\frac {343717911}{1863680}}\,{n}^{9}-{\frac {265635477}{116480}}\,{n}^{
8}
$$
$$
+{\frac {93771627}{7280}}\,{n}^{7}-{\frac {45012561}{1040}}\,{n}^{6}+{\frac {43093563}
{455}}\,{n}^{5}-{\frac {77940021}{560}}\,{n}^{4}+{\frac {62196357}{455}}\,{n}^{3}-{
\frac {39213532}{455}}\,{n}^{2}+{\frac {14354176}{455}}\,n-{\frac {463661}{91}} )
\quad .
$$

Now let's standardize and take the limit as $n$ goes to infinity. The limiting kurtosis is
$$
\kappa_4=\lim_{n \rightarrow \infty} \frac{m_4(n)}{m_2(n)^2} = \frac{21}{5}=4.2 \quad,
$$
so now the Gepner statistics is even more  {\it leptokurtic}.

The standardized sixth moment tends to 
$$
\kappa_6=\lim_{n \rightarrow \infty} \frac{m_6(n)}{m_2(n)^3} = {\frac {279}{7}}  \quad .
$$

The standardized eighth moment tends to
$$
\kappa_8=\lim_{n \rightarrow \infty} \frac{m_8(n)}{m_2(n)^4} = {\frac {3429}{5}} \quad .
$$
Similarly, $\kappa_{10}=206955/11$ and $\kappa_{12}=343717911/455$.
A quick {\it google} search for  ``21/5, 279/7 , 3429/5'' revealed that these are the 
$4th$, $6th$, and $8th$ moments 
of  the standard {\bf Logistic distribution}, whose probability generating function is
$$
\frac{\pi}{4 \sqrt{3} } sech^2 (\frac{\pi x}{2 \sqrt{3}}) \quad ,
$$
and whose moments are $(2^n-2) \cdot |B_n| \cdot 3^{n/2}$, where $B_n$ are the Bernoulli numbers.

So the first $12$ moments perfectly agree, and I am sure that they all do,
and hence I can assert with great confidence the interesting

{\bf Fact}:
Doron Gepner's Statistics on the set of words $\W(n,n,n)$ has mean $\frac{n^3}{2}$, variance $\frac{n^4}{4}$ and
its limiting (scaled) distribution is the  {\bf logistic distribution}.

\vfill\eject

Let me conclude with

{\bf YOM HULEDET SHISHIM SAMEACH, DORON!} \quad .

{\bf Acknowledgment}

I wish to thank my beloved servant, Shalosh B. Ekhad, for its diligent computations,
and to Dominique Foata for helpful discussions  way back in 1987.

{\bf References}

[A] George Andrews, {\it ``The Theory of Partitions''}, Cambridge University Press, 1984.
Originally published, in 1976, by Addison-Wesley.

[BZ] Andrew Baxter and Doron Zeilberger,
{\it The Number of Inversions and the Major Index of Permutations are Asymptotically Joint-Independently Normal},
The Personal Journal of Shalosh B. Ekhad and Doron Zeilberger, Feb. 4, 2011,\hfill\break
{\tt http://www.math.rutgers.edu/\~{}zeilberg/mamarim/mamarimhtml/invmaj.html} \quad .

[CJZ] E. Rodney Canfield, Svante Janson, and Doron Zeilberger, {\it The Mahonian Probability Distribution on Words is Asymptotically Normal},
Advances in Applied Mathematics {\bf 46} (2011), 109-124.

[D] Persi Diaconis, {\it ``Group Representations in Probability and Statistics''}, Institute of Mathematical Statistics,
Hayward, CA, 1988.

[F] Dominique Foata, {\it On the Netto inversion number of a sequence}, Proceedings of the American Mathematical
Society {\bf 19} (1968), 236-240.

[G] Doron Gepner, {\it private communication}, 1987.

[MW] H.B. Mann and D.R. Whitney, {\it On a test whether one of two random variables is stochastically larger than other},
Annals of Mathematical Statistics {\bf 18} (1947), 50-60.

[N] Brian Nakamura, {\it ``Computational Methods in Permutation Patterns''}, PhD thesis, Rutgers University, \hfill\break
{\tt http://www.math.rutgers.edu/\~{}zeilberg/Theses/BrianNakamuraThesis.pdf} \quad .

[Z1] Doron Zeilberger,
{\it The Automatic Central Limit Theorems Generator (and Much More!) },
in:"Advances in Combinatorial Mathematics: Proceedings of the Waterloo Workshop in Computer Algebra 2008
in honor of Georgy P. Egorychev", chapter 8, pp. 165-174,
(I.Kotsireas, E.Zima, eds. Springer Verlag, 2009), \hfill\break
{\tt http://www.math.rutgers.edu/\~{}zeilberg/mamarim/mamarimhtml/georgy.html} \quad .

[Z2] Doron Zeilberger, {\it HISTABRUT:A Maple Package for  Symbol-Crunching in Probability theory},
The Personal Journal of Shalosh B. Ekhad and Doron Zeilberger, \hfill\break
{\tt http://www.math.rutgers.edu/\~{}zeilberg/mamarim/mamarimhtml/histabrut.html} \quad .

\bigskip
\bigskip
\hrule
\bigskip
Doron Zeilberger, Department of Mathematics, Rutgers University (New Brunswick), Hill Center-Busch Campus, 110 Frelinghuysen
Rd., Piscataway, NJ 08854-8019, USA. \hfill \break
zeilberg at math dot rutgers dot edu \quad ;  \quad {\tt http://www.math.rutgers.edu/\~{}zeilberg/} \quad .
\bigskip
\hrule
\bigskip
Exclusively published in The Personal Journal of Shalosh B. Ekhad and Doron Zeilberger  \hfill \break
({ \tt http://www.math.rutgers.edu/\~{}zeilberg/pj.html})
and {\tt arxiv.org} \quad . 
\bigskip
\hrule
\bigskip
{\bf March 31, 2016}

\end